\documentstyle[12pt]{article}

\newcommand{\sect}[1]{\setcounter{equation}{0}\section{#1}}

\newcommand{\eq}{\begin{equation}}
\newcommand{\eqa}{\begin{eqnarray}}
\newcommand{\en}{\end{equation}}
\newcommand{\ena}{\end{eqnarray}}
\newcommand{\enn}{\nonumber \end{equation}}

\def\sk{\vskip .4cm}
\def\noi{\noindent}

\def\al{\alpha}

\def\Cb{\bar{C}}

\def\epsi{\varepsilon}

\def\de{\delta}

\def\part{\partial}

\def\R#1#2{ R^{#1}_{~~~#2} }
\def\PA#1#2{ P^{#1}_{A~~#2} }

\def\Rb{{\mbox{\bf R}}}

\def\Rbold{\mbox{\scriptsize \bf R}}

\def\Rh{{\hat R}}

\def\Cb{\bf \mbox{\boldmath $C$}}

\def\T#1#2{ T^{#1}_{~~#2} }

\def\rminus{q^{-1}}

\def\D{\Delta}

\def\n2{{{N+1} \over 2}}
\def\ap{a^{\prime}}
\def\bp{b^{\prime}}

\def\Dc{{\cal D}}

\def\Ntwo{{N\over 2}}

\def\square{{\,\lower0.9pt\vbox{\hrule \hbox{\vrule height 0.2 cm
\hskip 0.2 cm \vrule height 0.2 cm}\hrule}\,}}

\def\QED{{\vskip -0.49cm\rightline{$\Box$}}}

\def\#{\sharp}
\def\*{\star}
\def\a*{*^\sharp}
\def\b*{{\*}^{\sharp}}
\def\x{\hspace{.2mm}\mbox{\boldmath $\scriptscriptstyle{}^{{}^{_{{}_{{}^{\!\times}}}}}$}}
\def\xx{{\scriptstyle{\times}}}
\def\cvd{{\vskip -0.49cm\rightline{$\Box\!\Box\!\Box$}}\sk}
\def\le{\langle}
\def\re{\rangle}

\def\g{\mbox{\sl g}}

\begin{document}

\begin{titlepage}
\rightline{LBNL-41692}
\rightline{April 1998}
\vskip 2em
\begin{center}{\large\bf  Real forms of quantum orthogonal groups,\\[.4em]
{\mbox{\boldmath $q$}}-Lorentz groups in any dimension}\\[5em]
Paolo Aschieri \\[.5em]
{\sl Theoretical Physics Group, Physics Division\\
Lawrence Berkeley National Laboratory, 1 Cyclotron Road \\
Berkeley, California 94720, USA.}\\[10em]
\end{center}

\begin{abstract}

We review known real forms of the quantum orthogonal groups $SO_q(N)$.
New $*$-conjugations are then introduced and we
contruct all real forms of quantum orthogonal groups. We thus give an $RTT$ 
formulation of the $*$-conjugations on $SO_q(N)$ that is complementary to the 
${\cal U}_q(\mbox{\sl g})$ 
$*$-structure classification of Twietmeyer \cite{Twietmeyer}.
In particular we easily find and describe the real forms 
$SO_q(N-1,1)$ for any value of $N$.
Quantum subspaces of the $q$-Minkowski space are analized.
\end{abstract}

\vskip 3.2cm

\noi MSC-class: 81R50; 17B37~~~~~~~~~~~~~~~~~~~~~~~~~~~
{}~~~~~~~~~~~~~~~~~~~~~~~~~~~~~~~~~~~~~~~~~
\vskip .2cm
\noi \hrule
\vskip .2cm

\noi{\small e-mail: aschieri@lbl.gov
}

\vskip.2cm


\end{titlepage}

\newpage
\setcounter{page}{1}

\sect{Introduction}
In the last years there has been an intense study of noncommutative 
deformations of
the Lorentz group. A deformed Lorentz group acts on a deformed noncommutative
Minkowski spacetime that has thus the same symmetry richness as in the 
classical case.    
An important aspect is that non commuting coordinates give rise to 
indetermination 
relations and discretization, this might be a realistic picture of 
how spacetime behaves 
at short distances (see for example  \cite{Wess}).

There are many deformations of the Lorentz group. In this paper we are 
concerned with
the standard FRT \cite{FRT} deformation of the orthogonal groups.
There, the complex orthogonal quantum groups $SO_q(N)$ are defined for 
any $N> 2$ 
and various real forms are studied, in particular $SO_q(2n,1)$ with 
$q\in {\bf\mbox{R}}$. 
A new real form $SO_q(n+1,n-1)$ with $|q|=1$  was then described in 
\cite{Firenze1}. 
We present here the real forms $SO_q(2n-1,1)$ with $q\in {\bf\mbox{R}}$ 
and therefore construct
$q$-Lorentz groups in any dimensions. More in general, for 
$q\in {\bf\mbox{R}}$, 
we can obtain any desired 
signature, $SO_q(l,m),~ l+m=N$ as well as deformations of the real form 
$SO^*(2n)$.
\sk
This letter is organized as follows. 
We first briefly recall the $R$-matrix construction of the  
orthogonal quantum groups, we then discuss in the undeformed case the
equivalence between real forms and $*$-structures (conjugations). It is the 
$*$-structure formulation that has a straighforward generalization to the 
quantum case, see however \cite{Dobrev} for a different 
approach.
We end the section recalling
some of the $*$-conjugations defined in \cite{FRT}. 
Section 3 contains the main results, we
find an involution $\sharp$ that is an $SO_q(N)$ automorphism 
(algebra and coalgebra map compatible with the antipode
of $SO_q(N)$). The composition $\a*=\# {\scriptstyle{{}^{{}_{\circ}}}}*$ of 
$\sharp$ with the $*$-conjugations defined in \cite{FRT}
gives  new star structures that lead to the real forms:  
$SO_q(2n-1,1)$ with $q\in {\bf\mbox{R}}$, 
$SO_q(n,n+1)$ with $|q|=1$ and $SO_q(n+1,n-1)$ with $|q|=1$.
The conjugations on the quantum orthogonal plane that are 
associated to these $\a*$-conjugations are also easily found.
This procedure can be reiterated with another set of involutive automorphisms   
$\natural$, $\# {\scriptstyle{{}^{{}_{\circ}}}}\natural$, that composed with 
$*$  
give, for $q\in {\bf\mbox{R}}$, 
pseudo-orthogonal $q$-groups
with any signature as well as the $q$-group{\sl s} $SO^*_q(2n)$. Comparison of 
our results with 
\cite{Twietmeyer}  shows that  via the involutive 
automorphisms $\sharp$, $\natural$ we obtain a complete classification of all 
the real forms of 
orthogonal $q$-groups.
In Note 2 and Note 1 we remark that this is the same 
construction as in the classical case, 
where all real forms of $SO(N,\mbox{\bf C})$ are classified via 
involutive automorphism $\#,\,\natural$ of the compact real form 
$SO(N, \mbox{\bf R})$.  
These involutive  automorphisms $\#\,(\natural)$ of 
$SO(N, \mbox{\bf R})$ can be realized
via  matrices $\Dc \in O(N,\mbox{\bf R})$:
\(
{}~\# : SO(N, \mbox{\bf R})\rightarrow SO(N, \mbox{\bf R}),
{}~T^{\#}=\Dc T\Dc^{-1}
\)
with $\Dc^2=\pm${\bf 1}.
Similarly in the quantum case the $SO_q(N,\Rb)$ involutive automorphisms 
${\sharp},\,\natural$
are realized via matrices (with complex entries) that satisfy the $RTT$ 
and orthogonality 
relations of $SO_q(N)$, that have $q$-determinant $=\pm 1$, that squared 
equal $\pm${\bf 1} 
and that moreover satisfy the 
reality conditions of $SO_q(N,\mbox{\bf R})$.
In Note 2 we also explicitly classify the involutive automorphisms of 
$SO_q(N,\Rb)~(\forall\,N\not=8)$.

We conclude by observing that even if the classical  inclusion $SO(N)\subset$
$SO(N+1)$ no more holds at the quantum level, one still has that the
$(N+1)$-dimensional $q$-orthogonal plane includes 
the $N$-dimensional $q$-orthogonal plane.  In particular the $q$-Minkowski 
spacetime contains the $q$-euclidean space.

\section{$SO_q(N)$ and the real forms $SO_{q}(N,\bf\mbox{R})$, 
$SO_{q}(n,n)$,  $SO_{q}(n+1,n)$}

The 
quantum group $SO_{q}(N)$ is freely generated by the non-commuting
matrix elements $\T{a}{b}$ (fundamental representation, $a,b= 1,\ldots N$) and
the unit element $I$, modulo  the 
quadratic $RTT$ and $CTT$ (othogonality) 
relations discussed below and the condition det$_qT=1$ where det$_q$ is the 
quantum determinant 
defined via the quantum epsilon tensor $\epsi_q^{i_1,...i_N}$ in \cite{Fiore}. 
The noncommutativity is controlled by the $R$  
matrix:
$\R{ab}{ef} \T{e}{c} \T{f}{d} = \T{b}{f} \T{a}{e} \R{ef}{cd}$ i.e.:
\eq
R_{12}T_1T_2=T_2T_1R_{12}~.
\label{RTT}
\en
The ortogonality relations are
$C^{bc} \T{a}{b}  \T{d}{c}= C^{ad} I$,
$C_{ac} \T{a}{b}  \T{c}{d}=C_{bd} I$
that we rewite in matrix notation as:
\eq
TCT^t=C\,I ~~,~~~~T^tCT=C\,I\label{Torthogonalitymat}~.
\en
Following \cite{FRT} we define prime indices  as $\ap \equiv N+1-a$.
The indices run over $N$ values, for $N$ even we write $N=2n$, otherwise 
$N=2n+1$, in this case
we define $n_2  \equiv \n2$. Notice that for   $N=2n$ we have $n'=n+1$.
We also define a 
vector $\rho_a$  as:
\eq
(\rho_1,...\rho_N)=\left\{ \begin{array}{ll}
         (\Ntwo -1, \Ntwo -2,...,
{1\over 2},0,-{1\over 2},...,-\Ntwo+1)
                   & \mbox{for $SO(2n+1)$} \\
           (\Ntwo -1,\Ntwo -2,...,1,0,0,-1,...,-\Ntwo+1) & \mbox{for
$SO(2n)$}
                                             \end{array}
                                    \right.
\en
Then the (antidiagonal) metric is :
\eq
C_{ab}= q^{-\rho_a} \de_{a\bp} \label{metric}
\en
\noi its inverse $C^{ab}$
satisfies $C^{ab} C_{bc}=\de^a_c=C_{cb} C^{ba}$, it is easily seen that  
the matrix elements of the metric and the inverse metric coincide.
\sk
The nonzero complex components
of the $R$ matrix are (no sum on repeated indices):
\eqa
& &\R{aa}{aa}=q , ~~~~~~~~~~~~~~\mbox{\footnotesize
 $a \not= {n_2}$ } \nonumber\\
& &\R{a\ap}{a\ap}=q^{-1} ,  ~~~~~~~~~~\mbox{\footnotesize
 $a \not= {n_2}$ } \nonumber\\
& &\R{{n_2}{n_2}}{{n_2}{n_2}}= 1\nonumber\\
& &\R{ab}{ab}=1 ,~~~~~~~~~~~~\mbox{\footnotesize
 $a \not= b$,  $\ap \not= b$}\label{Rnonzerososp}\\
& &\R{ab}{ba}=q-q^{-1} , ~~~~~~~\mbox{\footnotesize
$a>b, \ap \not= b $}\nonumber\\
& &\R{a\ap}{\ap a}=(q-q^{-1})(1- q^{\rho_a-\rho_{\ap}})
=(q-q^{-1})[1-C^{a'a}C_{a'a}] ,
{}~~~\mbox{\footnotesize
$a>\ap $}\nonumber\\
& &\R{a\ap}{b \bp}=-(q-q^{-1})
q^{\rho_a-\rho_b}=-(q-q^{-1})C^{a'a}C_{bb'} , ~~~~~~~
\mbox{\footnotesize $~~a>b ,~ \ap \not= b $}\nonumber
\ena
The coproduct, counit and antipode are as usual defined by:
$\D(\T{a}{b})=\T{a}{c}\otimes\T{c}{b}$ $\epsi(\T{a}{b})=\delta^a_b$ 
and $\kappa(T)=CT^tC$
i.e.  $\kappa(T^a{}_b)=C^{ae}\T{f}{e}C_{fb}$.
\sk
Quantum orthogonal planes are defined by the
commutations
\eq
\PA{ab}{cd} x^c x^d=0\label{qplane}
\en
where the $q$-antisymmetrizer $P_A$ is given by
$
P_A={1 \over {q+\rminus}} [-\Rh+rI-(q- q^{1-N})P_0] $ and 
${P_0}^{ab}_{~cd}= (C_{ef} C^{ef})^{-1} C^{ab} C_{cd}$.
The symmetry of the  $q$-planes under the $q$-groups is expressed via
the coaction $\delta$: $\delta(x^a)=T^a{}_b\otimes x^b$.
\sk
\noi{\bf Star Structures and  Real Forms  }
\sk
The Hopf algebra $SO_q(N)$ 
can be interpreted as the deformation of the algebra
of functions on a group manifold only introducing a 
$*$-structure on $SO_q(N)$ (the analogue of complex conjugation). 
This procedure leads to the quantum groups
$SO_q(N, \Rb)=Fun_q(SO(N,\Rb))$, $SO_q(N-1,1)=Fun_q(SO(N-1,1))$  and 
more in general $SO_q(l,m)=Fun_q(SO(l,m))$ and $SO_q^*(2n)=Fun_q(SO^*(2n))$.
A $*$-structure or $*$-conjugation  on a Hopf algebra $A$ 
(and in particular on $SO_q(N)$) 
is an algebra anti-automorphism $(\eta ab)^*=\overline{\eta}\,b^*a^*~\forall 
a,b\in A,~ 
\forall \eta\in$ {\bf C};
coalgebra automorphism $\D{\scriptstyle{{}^{{}_{\circ}}}}*=
(*\otimes *) {\scriptstyle{{}^{{}_{\circ}}}}\D$, 
$\epsi {\scriptstyle{{}^{{}_{\circ}}}}*= \epsi$ 
and involution  $*^2=id$. It follows that
$*{\scriptstyle{{}^{{}_{\circ}}}}\kappa^{-1}=
\kappa{\scriptstyle{{}^{{}_{\circ}}}}*$
i.e. $[\kappa^{-1}(T)]^*=\kappa(T^{*})$.\footnote{{\sl Proof}: 
$~\de^a{}_cI=\epsi({T^{a}{}_{c}}^*)I=(\kappa\otimes id)\D({T^{a}{}_{c}}^*)
=\kappa({T^{a}{}_{b}}^*){T^{b}{}_{c}}^*$. Apply $*$ to get $\de^a{}_cI=
T^{b}{}_{c}[\kappa({T^{a}{}_{b}}^*)]^*$; multiply from the left with 
$\kappa^{-1}(T^{c}{}_{e})$ to get the thesis. \QED}
Two star structures $*$ and $*'$ 
are equivalent if there exist a Hopf algebra automorphism
$\alpha$ such that 
$*'=\al {\;\scriptstyle{{}^{{}_{\circ}}}}\!*{\scriptstyle{{}^{{}_{\circ}}}}
\al^{-1}$..

Let us clarify the relation 
between real forms of groups
and $*$-structures on Hopf algebras. 

Consider a  real form $\mbox{\sl g}_{\Rbold}$ of a complex Lie algebra  
$\mbox{\sl g}$ i.e. $\mbox{\sl g}=\mbox{\sl g}_{\Rbold}
\oplus i\,\mbox{\sl g}_{\Rbold}$. We exponentiate the elements of 
$\mbox{\sl g}_{\Rbold}$ to get the real form $G_{\Rbold}$ 
of the complex Lie group $G$. 
There is a one-to-one correspondence between real forms of the complex 
Lie algebra
$\mbox{\sl g}$ and star structures on the universal enveloping algebra 
${\cal U}(\mbox{\sl g})$. The $*$-operation that acts as 
minus the identity on $\mbox{\sl g}_{\Rbold}$ satisfies $[\chi,\chi']^*=
[{\chi'}^*,\chi^*]\;$ $\forall \chi,\chi'\in\mbox{\sl g}_{\Rbold}$
and is uniquely extended as an anti-linear, anti-multiplicative 
and involutive map on the Hopf algebra 
${\cal U}(\mbox{\sl g})$. We also have a corresponding $*$-structure on  
${\cal{F}}(G)$, the algebra of analytic functions on the complex group $G$.
We define
$* \,:\; f\mapsto f^*$  such that $f^*(g)\equiv 
\overline{f(g)} ~\forall g\in G_{\Rbold}$, then
$f^*$ is extended analytically to all $G$. Viceversa 
a $*$-structure on ${\cal{F}}(G)$ determines the following real form 
$G_{\Rbold}$ 
of $G$: $G_{\Rbold}=\{g\in G \;| ~f^*(g)=\overline{f(g)}\}$.
In other words, the algebra of regular functions on
 $G_{\Rbold}$ [for ex.  $SO(N,\Rb)$] is isomorphic to the algebra of regular 
functions on  $G_{\Rbold '}$ [for ex.  $SO(N-1,1)$] 
and to the algebra of analytic functions on the 
complex manifold $G$ [$SO(N,\Cb)$]; 
indeed any regular function  $f\in Fun(G_{\Rbold})$ can be 
analytically continued in a unique 
function $\hat f\in {\cal{F}}(G)$, then the restriction 
of $\hat{f}$ to the $G_{\Rbold '}$ sub-manifold of $G$
belongs to $Fun(G_{\Rbold '})$.
Only considering a $*$-structure we can understand if
the orthogonal matrix entries  $T^a{}_b$ 
generate functions on $SO(N,\Rb)$ or on $SO(N-1,1)$.

We have seen that a real form $\mbox{\sl g}_{\Rbold}$  of $\mbox{\sl g}$
is equivalent to a $*$-structure on ${\cal U}(\mbox{\sl g})$ and to a 
$*$-structure on ${\cal{F}}(G)$, the algebra of analytic functions on $G$.
The $*$-structures on ${\cal{F}}(G)$ and ${\cal U} (\mbox{\sl g})$ can be directly related
using the duality $\le~,\re$ between these two Hopf algebrae.
Since
$f^*(g)=\overline{f(g)}$ $\forall g\in G_{\Rbold}$ we have, write $g=e^{\chi}$,
$\le\chi,f^*\re=\overline{\le\chi , f\re}~\forall 
\chi\in \mbox{\sl g}_{\Rbold}$. 
More in general
$\;\forall\psi\in {\cal{U}}(\mbox{\sl g})\,,~
\forall f\in {\cal{F}}(G)$,
\eq
\le\psi,f^*\re={\overline{\le\kappa(\psi)^{\,*}{},f\re}}
\label{U*A}~.
\en
{\sl Proof}:$~$ write 
\eq
\le\psi,f^*\re={\overline{\le {\cal{O}}(\psi){},f\re}}\label{defO}
\en
we want to find the map ${\cal{O}}$ implicitly defined by the above 
relation and 
such that ${\cal{O}}$ on $\mbox{\sl g}_{\Rbold}$ is the identity. 
From (\ref{defO}) it is easy to see that ${\cal{O}}$ is 
antilinear and multiplicative, this uniquely defines 
$\cal{O}=*{\scriptstyle{{}^{{}_{\circ}}}}\kappa$, 
indeed  
$* {\scriptstyle{{}^{{}_{\circ}}}}\kappa(\chi)=-\chi^*=\chi$.
\QED

Finally it is straighforward to see that if $*$ is equivalent to $*'$ 
i.e., if there exists
an automorphism $\alpha$ of ${\cal{U}}(\mbox{\sl g})$ [or ${\cal{F}}(G)$]
such that 
$*'=\al {\;\scriptstyle{{}^{{}_{\circ}}}}\!*{\scriptstyle{{}^{{}_{\circ}}}}
\al^{-1}$,
then $\mbox{\sl g}_{\Rbold}$ is isomorphic to  $\mbox{\sl g}_{\Rbold'}$ and 
$G_{\Rbold}$ is isomorphic to $G_{\Rbold'}$; the viceversa holds as well.

\sk
In the quantum case, following \cite{FRT}, on orthogonal $q$-groups 
a conjugation can be defined
\sk
$\bullet$~~ trivially as $T{\x}=T$ i.e. $\T{a}{b}{\!\x}=\T{a}{b}$. 
Compatiblility
with the
$RTT$ relations (\ref{RTT}) requires 
${\bar R}_{q}=R^{-1}_{q}=
R_{q^{-1}}$,
i.e. $|q|=1$. Then the $CTT$ relations are invariant under
$*$-conjugation.  The corresponding real forms are
$SO_q(n,n;\bf{\mbox{R}})$, $SO_q(n+1,n;\bf{\mbox{R}})$. To prove this,
recall that
in the classical limit 
a $*$ structure on $G$ determines the  real form  
$G_{\Rbold}=\{g\in G \;| ~f^*(g)=\overline{f(g)}~ \forall f\in{\cal{F}}(G)\}$;
therefore in the $q\rightarrow 1$ limit the conjugation $*$ becomes
the usual complex conjugation 
${}^{{}^{\overline{~~}}}$ on functions on $G_{\Rbold}$. 
In our case  the condition $T{\x}\!(g)=\overline{T(g)}$ reads 
$T(g)=\overline{T(g)}$. We have real matrices $T(g)$ 
orthogonal with respect to the metric $C_{q=\!1}^{ab}=\de^{ab'}$.
In the new basis $T'=MT$, where $M^a{}_b=
{1\over\sqrt{2}} (\de^a{}_b+\epsilon^a\de^a{}_{b'})$
and $\epsilon^a=1$ if $a<a'$, $\epsilon^a=-1$ if $a>a'$,
$\epsilon^a=\sqrt{2}-1$ if $a=a'$, we have 
${T'}^tC'T'=C'I$ where $C'=(M^{-1})^tCM^{-1}$ is diagonal with 
signature $(n,n)$ 
or $(n+1,n)$..

\sk
$\bullet$~~ via the metric as
$T^{\*}=[\kappa(T)]^t$ i.e. $T^{\*}=C^tTC^t$.
The condition on $R$ is
${\bar \R{ab}{cd}}=\R{dc}{ba}$, which
happens for $q \in $ {\bf\mbox{R}}.
Again the $CTT$ relations are ${\*}$-invariant.
The metric on a real basis (we will see the explicit construction 
later on) has in the $q\rightarrow 1$ limit
compact signature $(+,+,...+)$ so that the real form is $SO_{q}(N;\Rb)$.
\sk

The quantum orthogonal  group co-acts on the quantum orthogonal plane, 
and may induce an associated 
conjugation on the $q$-plane as well.  
More precisely  a conjugation on the quantum orthogonal plane ---i.e. 
an anti-linear 
anti-multiplicative involution of the $q$-plane algebra--- is compatible 
with a conjugation on its $q$-symmetry group if  the coaction
$\delta$  of the $q$-group on the $q$-plane: $\delta(x^a)=
T^a{}_b\otimes x^b$ satisfies $\delta(x^*)=T^*\otimes x^*\equiv\delta^*(x)$.
The above two $q$-groups conjugations have  unique
(up to an overall phase) associated  $q$-plane conjugations,
they are respectively  $(x^a){\x}=x^a$ and $(x^a)^{\*}=C_{ba}x^b$.

\section{Real forms $SO_{q}(N-1,1)$, $SO_q(l,m)$ and $SO_q^*(2n)$ }

We first notice that if we have an involution  
${\sharp}$ that is a 
Hopf algebra automorphism (algebra and 
coalgebra map compatible with the antipode: 
$\kappa(a^{\sharp})=[\kappa(a)]^{\sharp}\,{}$) and that commutes with a 
conjugation $*$, 
then the composition of these two involutions: 
${*^{\sharp}}\equiv {\sharp}{\scriptstyle{{}^{{}_{\circ}}}}*= 
*{\scriptstyle{{}^{{}_{\circ}}}}{\sharp}$ is 
again a conjugation. [Hint: it is 
trivially antilinear, antimultiplicative and compatible with the coproduct:
$\D(T^{*^{\sharp}})=T^{*^{\sharp}}\otimes T^{*^{\sharp}}$.
It is involutive because 
$({\sharp}{\scriptstyle{{}^{{}_{\circ}}}}*) {\scriptstyle{{}^{{}_{\circ}}}}
(\sharp{\scriptstyle{{}^{{}_{\circ}}}}{*})=id \Leftrightarrow 
{\sharp}{\scriptstyle{{}^{{}_{\circ}}}}*= 
*{\scriptstyle{{}^{{}_{\circ}}}}{\sharp}$].
We now find an involution $\sharp$ that comutes with
$\xx$ and $\*$ as defined above.
\sk
Define the map ${\sharp}$ on the generators of $SO_q(N)$ as:
\eq
T^{\sharp}=\Dc T{\Dc}^{-1}~~~\mbox{i.e.}~~~
(T^a{}_b)^{\sharp}=\Dc^a{}_eT^e{}_f{\Dc^f{}_b}^{-1}\label{Def}
\en
and extend it by linearity and multiplicativity to all $SO_q(N)$.
The entries of the $N$-dimensional $\Dc$  matrix are
\eq
\begin{array}{cc}
\Dc=\scriptsize{
\left(  \begin{array}{cccccccc}
{1}&{}&{}&{}&{}&{}&{}&{}\\
{}&{\cdot\cdot\cdot}&{}&{}&{}&{}&{}&{}\\
{}&{}&{1}&{}&{}&{}&{}&{}\\
{}&{}&{}&{0}&{1}&{}&{}&{}\\
{}&{}&{}&{1}&{0}&{}&{}&{}\\
{}&{}&{}&{}&{}&{1}&{}&{}\\
{}&{}&{}&{}&{}&{}&{\cdot\cdot\cdot}&{}\\
{}&{}&{}&{}&{}&{}&{}&{1}
\end{array}\right)}~,~~
&
\Dc=\scriptsize{
\left(  \begin{array}{ccccccc}
{1}&{}&{}&{}&{}&{}&{}\\
{}&{\cdot\cdot\cdot}&{}&{}&{}&{}&{}\\
{}&{}&{1}&{}&{}&{}&{}\\
{}&{}&{}&{-1}&{}&{}&{}\\
{}&{}&{}&{}&{1}&{}&{}\\
{}&{}&{}&{}&{}&{\cdot\cdot\cdot}&{}\\
{}&{}&{}&{}&{}&{}&{1}
\end{array}
\right)}\\
\mbox{\small for $N$ even} &
\mbox{\small for $N$ odd}
\end{array}\label{Dcmat}
\en 
In the $N=2n$ case the $\Dc$ matrix exchanges the index $n$ with the 
index $n+1$, in the $N=2n+1$ case $\Dc$ change the sign of the entries  of 
the $T$ matrix as many times as the index $n_2={N+1\over 2}$ appears.
Since $\Dc^2={\bf 1}$ we immediately see that $\sharp$ is an involution.

We now prove that $\sharp$  is compatible with the algebra structure, 
i.e. it is compatible with the $RTT$ and $CTT$ relations.
Use (\ref{Rnonzerososp}) to prove that
\eq
\Dc^{-1}_1\Dc^{-1}_2R_{12}\,\Dc_1\Dc_2=R_{12} ~~\mbox{ i.e. }~~ 
R_{12}\,\Dc_1\Dc_2=\Dc_2\Dc_1R_{12}.{\label{DDRDD}}
\en
We also have 
\eq
\Dc^tC\Dc=C ~~~,~~~~~ ~\Dc C\Dc^t=C\label{DCD}
\en
(the second equation in (\ref{DCD}) follows from the first). Therefore $\Dc$ 
is a matrix 
of complex entries
that satisfies the $RTT$ and orthogonality relations. It is then
trivial to show that $\#$ is compatible with the algebra structure,
indeed $T^{\#}=\Dc T\Dc^{-1}$ is a product of $q$-orthogonal matrices whose
matrix entries mutually commute, therefore $T^{\#}$ itself is a $q$-orthogonal matrix. 
We conclude that $R_{12}T^{\#}_1 T^{\#}_2=T^{\#}_2 T^{\#}_1 R_{12}$,
$T^{\#}C T^{\#\,t}=C\,I$, $T^{\#\,t}C T^{\#}=C\,I$ are equivalent to (\ref{RTT})
and (\ref{Torthogonalitymat}) and therefore $\sharp$ is an algebra automorphism of $SO_q(N)$.
An explicit proof is also instructive, for example we have: 
$(TC T^t)^{\sharp}=\Dc T\Dc\, C\, \Dc^t T^t \Dc^t=
\Dc (T C T^t) \Dc^t=\Dc C\Dc^t\,I=C\,I$.

The compatibility of  $\sharp$ with the  coalgebra
structure is easily checked:
 $\D(T^{\sharp})=\Dc\D(T)\Dc^{-1}=T^{\sharp}\otimes T^{\sharp}$,
$\epsi(T^{\#})=\epsi(T)$. Compatibility with the antipode: 
$\kappa(a^{\sharp})=[\kappa(a)]^{\sharp}\,{}$ follows from compatibility with
the algebra and coalgebra structure (cf. footnote {\small 1}).
We now show that the two conjugations defined in the previous  
section commute with ${\sharp}$. For the second conjugation, defined by
$T^{\*}=[\kappa(T)]^t=C^tTC^t$,  we have
$(T^{\sharp})^{\*}=(T^{\*})^{\sharp} ~\Leftrightarrow $ $\overline{\Dc}C^tTC^t
\overline{\Dc}^{-1}=
C^t\Dc T\Dc^{-1}C^t$. 
For an arbitrary $\Dc$ matrix this last relation is equivalent to 
(use Schur lemma)  
$\overline{\Dc}C^t=\mbox{\sl const}\cdot C^t\Dc$ i.e. 
\eq
\overline{\Dc}=\mbox{\sl const}\cdot C^t\Dc C^t\label{uffi}~.
\en
In our case $\Dc$ is given by (\ref{Dcmat}) and satisfies 
$\overline{\Dc}=C^t\Dc C^t$ so that $(T^{\sharp})^{\*}=(T^{\*})^{\sharp}$.
The two maps $\*$ and $\sharp$
not only commute when applied to the $T^a{}_b$ 
matrix entries, they also commute when applied 
to  any element of the $q$-group because they are respectively
multiplicative  and antimultiplicative. 
{}For the first conjugation, defined by $T{\x}=T$,
the proof that ${\sharp}{\scriptstyle{{}^{{}_{\circ}}}}\xx=$ 
$\xx{{\scriptstyle{{}^{{}_{\circ}}}}}$${\sharp}$ is straighforward.
\sk
Associated to  ${{\xx}^{\sharp}}$ and  ${{\*}^{\sharp}}$ we
can construct  the conjugations that act on the quantum orthogonal plane 
and are compatible with the coaction  $x^a\rightarrow T^a{}_b\otimes x^b$. 
These conjugations respectively are: 
$(x^a){{\x}^{{}^{\!\!\!\sharp}}}=\Dc^a{}_bx^b$ and  
$(x^a)^{{\*}^{\sharp}}=C_{ba}\Dc^b{}_e x^e$. 
Notice that these $q$-plane maps are well defined conjugations. 
This is so because $\#$, defined by $x^{\#}\equiv \Dc x$, is an involution: 
$\Dc^2=\mbox{\bf 1}$, and also on the $q$-plane $\#$ commutes with $\xx$ or $\*$.
\sk
The above results are particular examples of the following theorem:
\sk
\noi{\bf Theorem 1}  $~$ Consider a  map ${\sharp}$ on $SO_q(N)$ 
defined by (\ref{Def}) where $\Dc$ is an arbitrary matrix satisfying 
(\ref{DDRDD}),
(\ref{DCD}) and $\Dc^2=\pm${\bf{1}}. 
Then $\#$ is an involutive automorphism of $SO_q(N)$.  Moreover the 
compositions 
${\xx}^{\sharp}$
and ${{\*}^{\sharp}}$ of the conjugations $\xx$ and $\*$ with the 
automorphism ${\sharp}$ is 
again a conjugation if and only if, respectively, 
${\sharp}{\scriptstyle{{}^{{}_{\circ}}}}\xx 
=\xx{{\scriptstyle{{}^{{}_{\circ}}}}}$${\sharp}$
and 
${\sharp}{\scriptstyle{{}^{{}_{\circ}}}}
{\*}={\*} {\scriptstyle{{}^{{}_{\circ}}}}{\sharp}$. 
These commutations hold if and only if
\eq
\left.
\begin {array}{c}
\overline{\Dc}=\Dc~~~\mbox{ with } \Dc^2=\mbox{\bf 1}\nonumber\\
\overline{\Dc}=-\Dc~~~\mbox{ with } \Dc^2=-\mbox{\bf 1}\nonumber\\
\end{array}\right\}~~\mbox{for } \xx\label{a0}~,
\en
\eq 
~~~~~\overline{\Dc}= C^t\Dc C^t\label{a}~~~~~~~~~~~~~~~~~~\mbox{ for } \*~.
\en
If $\Dc^2=${\bf{1}} 
then the associated maps  ${{\xx}^{\sharp}}$ or ${{\*}^{\sharp}}$ that act 
on the $q$-plane as:
\eq
x^{{\*}^{\sharp}}=C^t\Dc x ~~~,~~~~~x{{\x}^{{}^{\!\!\!\sharp}}}=\Dc x 
\label{plane}
\en
are well defined conjugations and are compatible with the $q$-group coaction.
\sk
\noi{\sl Proof}: $~$ 
The definition (\ref{Def}) and relations (\ref{DDRDD}), (\ref{DCD}) and 
$\Dc^2=\pm${\bf{1}}
assure that $\#$ is an automorphism and an involution of $SO_q(N)$. 
To prove that 
${\sharp}{\scriptstyle{{}^{{}_{\circ}}}}\* 
=\*{{\scriptstyle{{}^{{}_{\circ}}}}}$${\sharp}$ is equivalent to 
(\ref{a}), square relation  (\ref{uffi}) 
and use $\Dc^2=\overline{\Dc}^2=\pm${\bf{1}}
to deduce that {\sl const}$\;=\pm 1$. We cannot have 
{\sl const}$\;=-1$ because $\overline{\Dc}=-C^t\Dc C^t$ and $\Dc^tC\Dc=C$ 
would imply
$\overline{\Dc}=-{\Dc^t}^{-1}$ i.e. $\overline{\Dc}\Dc^t=-\mbox{\bf 1}$ that 
is absurd 
($\sum_j \overline{\Dc}_{ij}\Dc_{ij}\not= -1$).
Similarly $\overline{\Dc}=\pm\Dc$
is equivalent to 
${\sharp}{\scriptstyle{{}^{{}_{\circ}}}}\xx 
=\xx{{\scriptstyle{{}^{{}_{\circ}}}}}$${\sharp}$ 
[just consider {\bf{1}} instead of $C^t$ in (\ref{uffi})]. Since 
$\overline{\Dc}\Dc=
\mbox{\bf 1}$, 
$\overline{\Dc}=\Dc$ is only compatible with
$\Dc^2=\mbox{\bf 1}$, while 
$\overline{\Dc}=-\Dc$ is only compatible with
$\Dc^2=-\mbox{\bf 1}$.
Finally the $q$-plane map defined by $x^{\#}=\Dc x$ is an 
involutive automorphism of the quantum plane algebra if and only if 
$\Dc^2=\mbox{\bf 1}$. 
In this case 
$\overline{\Dc}=\Dc$ (for $\xx$) or  $C^t\overline{\Dc}=\Dc C^t$ (for $\*$) 
so that
$\#$ commutes with $\xx$ or $\*$  on the  $q$-plane. Compatibility of $\xx$ 
and $\*$ 
with the quantum group coaction $\delta$:
$\delta(x{\x})=T{\x}\otimes x{\x}$ and $\delta(x^{\*})=T^{\*}\otimes x^{\*}$ 
is straighforward.

\cvd

\noi{\bf Note 1}$~$ We can interpret the $\Dc$ matrix of Theorem 1 as a 
special $T$ matrix 
with complex entries. Then  conditions (\ref{DDRDD}),
(\ref{DCD}) and $\Dc^2=\pm${\bf{1}} simply state that $\Dc$ is a quantum 
orthogonal matrix,
moreover the constraints (\ref{a}) respectively state that $\Dc$ satisfies 
(modulo a minus sign in the $\xx$ case) the reality conditions
$T{\x}=T$, of $SO_q(n,n)$, $SO_q(n+1,n)$, and $T^{\*}=C^tTC^t$ of 
$SO_q(N,\Rb)$. 
\sk
In order to  study the real forms related to ${{\xx}}^{\sharp}$ and  
${{\*}^{\sharp}}$, where $\#$ is given by (\ref{Def}), ({\ref{Dcmat}), we  
find a 
linear transformation $x\rightarrow x'=Mx$, $T\rightarrow T'=MT M^{-1}$
such that the new coordinates $x'$ and $T'$ are real and the new metric 
$C'=(M^{-1})^tCM^{-1}$, 
at least in the $q\rightarrow 1$ limit, becomes 
diagonal.\footnote{{}For $q\not=1$ it is not possible to have 
a diagonal metric, indeed it is not possible 
to diagonalize $C$ via the transformation  $C'=(M^{-1})^tCM^{-1}$ because 
$C^t\not= C\Rightarrow C'^t\not=C'$. Of course in the $q\rightarrow 1$ limit  
$C^t= C$.}
Then the signature of the metric identifies the real form associated to the 
given 
$*$-conjugation. 
\sk
The transformation matrices $M$ and the corresponding transformed metrics 
in the $q=1$ limit for $\*$ and $\*^{\#}$ explicitly read:
\[
\begin{array}{clc}
\mbox{\small $\!\!\!\!\!\!\!\!N$ even }~~~~~~~ 
&~~~~~~~~~~~~~~,~~~~~~~~~~~~~~& ~~~
\mbox{ \small $N$ odd}
\end{array}
\nonumber
\]
\[
\left.
\begin{array}{lcl}
\!\!\scriptsize{
M=
{1\over \sqrt{2}}\left(  \begin{array}{cccccccc}
{1}&{}&{}&{}&{}&{}&{}&{\!1}\\
{}&{\!\cdot\cdot\cdot}&{}&{}&{}&{}&{\!\cdot\cdot\cdot}&{}\\
{}&{}&{\!1}&{}&{}&{\!1}&{}&{}\\
{}&{}&{}&{\!1}&{\!1}&{}&{}&{}\\
{}&{}&{}&{\!i}&{\!-i}&{}&{}&{}\\
{}&{}&{\!i}&{}&{}&{\!-i}&{}&{}\\
{}&{\!\cdot\cdot\cdot}&{}&{}&{}&{}&{\!\cdot\cdot\cdot}&{}\\
{i}&{}&{}&{}&{}&{}&{}&{\!-i}
\end{array}\right) }&,
&
\scriptsize{
M={1\over{\sqrt{2}}}
\left(  \begin{array}{ccccccc}
{1}&{}&{}&{}&{}&{}&{\!1}\\
{}&{\!\cdot\cdot\cdot}&{}&{}&{}&{\!\cdot\cdot\cdot}&{}\\
{}&{}&{\!1}&{}&{\!1}&{}&{}\\
{}&{}&{}&{\!\sqrt{2}}&{}&{}&{}\\
{}&{}&{\!i}&{}&{\!-i}&{}&{}\\
{}&{\!\cdot\cdot\cdot}&{}&{}&{}&{\!\cdot\cdot\cdot}&{}\\
{i}&{}&{}&{}&{}&{}&{\!-i}
\end{array}\right) }\nonumber
\\
\!\!C'=diag(1,... 1,1,1,1,... 1) &, & C'=diag(1,... 1,1,1,... 1) 
\nonumber\end{array}
~\right\}
{\begin{array}{c}
\mbox{\small $\!\!$for}\nonumber\\[-.4em]
\mbox{$\!\!\!\*$}
\end{array}}
\nonumber
\]
\vskip.001cm
\[
\left. 
\begin{array}{lcl}
\!\!\scriptsize{
M=
{1\over \sqrt{2}}\left(  \begin{array}{cccccccc}
{1}&{}&{}&{}&{}&{}&{}&{\!1}\\
{}&{\!\cdot\cdot\cdot}&{}&{}&{}&{}&{\!\cdot\cdot\cdot}&{}\\
{}&{}&{\!1}&{}&{}&{\!1}&{}&{}\\
{}&{}&{}&{\!\!1}&{\!1}&{}&{}&{}\\
{}&{}&{}&{\!\!-1}&{\!1}&{}&{}&{}\\
{}&{}&{\!i}&{}&{}&{\!-1}&{}&{}\\
{}&{\!\cdot\cdot\cdot}&{}&{}&{}&{}&{\!\cdot\cdot\cdot}&{}\\
{i}&{}&{}&{}&{}&{}&{}&{\!-i}
\end{array}\right) }&, 
&
\scriptsize{
M={1\over{\sqrt{2}}}
\left(  \begin{array}{ccccccc}
{1}&{}&{}&{}&{}&{}&{\!1}\\
{}&{\!\cdot\cdot\cdot}&{}&{}&{}&{\!\cdot\cdot\cdot}&{}\\
{}&{}&{\!1}&{}&{\!\!\!1}&{}&{}\\
{}&{}&{}&{\!\!i\sqrt{2}}&{}&{}&{}\\
{}&{}&{\!i}&{}&{\!\!\!-i}&{}&{}\\
{}&{\!\cdot\cdot\cdot}&{}&{}&{}&{\!\cdot\cdot\cdot}&{}\\
{i}&{}&{}&{}&{}&{}&{\!-i}
\end{array}\right) }\nonumber
\\
\!\!C'=diag(1,... 1,1,-1,1,... 1) &, & C'=diag(1,... 1,-1,1,... 1) 
\nonumber\end{array}
~\right\}
{\begin{array}{c}
\mbox{\small $\!\!$for}\nonumber\\[-.4em]
\mbox{$\!\*^{\#}$}
\end{array}}
\nonumber
\]
Therefore $\*$ gives the real form $SO_q(N,\Rb)$, while  
 $\*^{\#}$ gives the real form $SO_q(N-1,1)$.
The matrices corresponding to the ${\xx}^{\sharp}$ 
conjugation have a similar structure with nonzero
entries only in the diagonals.
For $N$ even we obtain the real form  $SO_q(n+1,n-1)$ with $|q|=1$. 
For $N$ odd we obtain the real form $SO_q(n,n+1)$ with $|q|=1$;
this is isomorphic to the real form $SO_q(n+1,n)$. 

\noi{\sl Proof}: We construct an automorphism $\al$ of $SO_q(N)$ such that 
$\al{\scriptstyle{{}^{{}_{\circ}}}}\xx^{\#}=
\xx{\scriptstyle{{}^{{}_{\circ}}}}${}$\al$.  
Consider 
\eq
\al\;:~SO_q(N)\rightarrow SO_q(N)\;,~~\al(T)\equiv ATA^{-1}\;,~~~A\equiv 
diag(i,...i,1,-i,...-i)
\en
The matrix $A$ (like $\Dc$) satisfies (\ref{DDRDD}) and (\ref{DCD}) so that 
$R\al(T_1)\al(T_2)=\al(T_2)\al(T_1)R$, $\al(T)^tC\al(T)=C$; therefore $\al$
is an automorphism of $SO_q(N)$.
We also have $\Dc A=-\overline{A}$ so that 
$\al({T{{\x}^{\!\!\#}}})=[\al(T)]{\x}$ indeed:
\eq
\al({T{{\x}^{\!\!\#}}})=\Dc\al(T)\Dc^{-1}=\Dc ATA^{-1}\Dc^{-1}=\overline{A}T
\overline{A}^{-1}=(ATA^{-1}){\x}=[\al(T)]{\x}\label{equival}~.
\en

\QED
\sk
We now find  real forms $SO_q(l,m)$ with $l$ and $m$ arbitrary by considering
other examples of $\Dc$ matrices that, as in Theorem 1, satisfy (\ref{DDRDD}),
(\ref{DCD}), $\Dc^2=${\bf{1}} and $\overline{\Dc}=\Dc$ (for ${\xx}$) or  
$\overline{\Dc}=C^tDC^t$ (for $\*$).
Following \cite{FRT} we first study the  matrices
$\Dc'={{diag}(\epsi_1,\ldots ,\epsi_N)}$, where $\epsi_j^2=1$,  
$\epsi_{j'}=\epsi_j$,
$j=1,\ldots ,N$ and $\epsi_{n+1}=1$ for $N=2n+1$ or $\epsi_{n}=\epsi_{n+1}=1$ 
for $N=2n$.
[In the even case $\epsi_j=-1\;\forall\, j<n$ gives 
$-diag(1,...1,-\epsi_n,-\epsi_{n+1},1,...1)$ 
that differs from
minus the identity only if 
$\epsi_{n}=\epsi_{n+1}=1$, this is why we impose $\epsi_n=\epsi_{n+1}=1$.
Similarly in the odd case, where moreover 
$\epsi_j=-1\;\forall\, j<n+1$ gives 
minus the $\Dc$ matrix in (\ref{Dcmat})].
The involutive automorphisms $\natural$ associated to the $\Dc'$ matrices are 
given by
$T^{\natural}=\Dc' T {\Dc'}^{-1}$; the corresponding conjugations are 
${\xx}^{\natural}$
and $\*^{\natural}$. Since the $\Dc'$ matrices commute also with
the $\Dc$ matrices in (\ref{Dcmat}) we have $(\Dc\Dc')^2=\mbox{\bf 1}$, 
$\#^\natural$ is an
involutive automorphism,  
${\xx}{\scriptstyle{{}^{{}_{\circ}}}}{\#^\natural}=
\#^\natural{\scriptstyle{{}^{{}_{\circ}}}}{\xx}$,
$\*{\scriptstyle{{}^{{}_{\circ}}}}\#^\natural=
\#^\natural{\scriptstyle{{}^{{}_{\circ}}}}\*$
and the compositions ${{\xx}^{\#}}^{\natural}$, ${\*^{\#}}^{\natural}$ are 
again conjugations. 
This holds both in the quantum group and in quantum plane case.

The conjugations ${{\xx}^{\natural}}$ and ${{{\xx}^{\#}}^{\natural}}$ 
still give the real forms
$SO_q(n+1,n)$, $SO_q(n,n)$, $SO_q(n+1,n-1)$ with $|q|=1$. 

\noi {\sl Proof}: rewrite (\ref{equival}) for the authomorphism $\al$ given by
$A={diag}(\sigma(\epsi_1),\ldots ,\sigma(\epsi_N))$ where $\sigma(\epsi_j)=1$ 
if $\epsi_j=1$,
$\sigma(\epsi_j)=i$ if  $\epsi_j=-1$ and $j<j'$,  $\sigma(\epsi_j)=-i$ if  
$\epsi_j=-1$ 
and $j>j'$.\QED 

The conjugations ${{\*}^{\natural}}$, studied in \cite{FRT},
give the real forms $SO_q(l,m)$ where, if $l+m$ is odd,  $l$ and $m$ 
are arbitrary numbers satisfying the only constraint $l+m=2n+1$, if $l+m$ is 
even $l$ and $m$ 
are both even.
Finally  real forms $SO_q(l,m)$ with $l+m=2n$ and $l$ and $m$ both odd are 
given by
the conjugations ${{\*^{\#}}^{\natural}}$.

In the $SO_q(2n+1)$ case there are $2^n$ matrices $\Dc'$  and  
we conclude that the number of real forms 
$SO_q(2n+1, \Dc')$ of $SO_q(2n+1)$ is $2^n$ for $q\in\Rb$. For $|q|=1$ we 
have just the real form $SO_q(n+1,n)$.
In \cite{Twietmeyer} 
Twietmeyer classifies real forms of ${\cal{U}}_q({\mbox{\sl g}})$ where 
${\mbox{\sl g}}$ is semisimple
(and more in general symmetrizable Kac-Moody); there are $2^n$ real forms of 
${\cal{U}}_q({so(2n+1)})$ for $q\in \Rb$ and just one real form for $|q|=1$.
Comparison with \cite{Twietmeyer} shows that  via  $\xx$, $\*$ and the $\Dc'$ 
matrices 
we have described {all}  real forms of $SO_q(2n+1)$.

One proceed similarly in the $SO_q(2n)$ case. We have $2^{n-1}$  matrices 
$\Dc'$ that give
$2^{n-1}$ $\;\*^{\natural}$-conjugations. The   conjugations 
$\*^{\#^{\natural}}$
are only $2^{n-2}$ since the matrices
$\Dc'_1={{diag}(\epsi_1,...\epsi_{n-1},1,1,\epsi_{n+2},...\epsi_{2n})}$ and
$\Dc'_2={{diag}(-\epsi_1,...-\epsi_{n-1},1,1,-\epsi_{n+2},...-\epsi_{2n})}$ 
give equivalent 
conjugations  $\*^{\#^{\natural_1}}$  and  
$\*^{\#^{\natural_2}}$.\\{\sl Proof}:$~$ 
The map $\al(T)=ATA^{-1}$ given by 
$A= diag(-1,...-1,1,...1)$ is an automorphism
of $SO_q(2n)$.  We also have
$C^t\Dc\Dc'_1A= \overline{A}C^t\Dc\Dc'_2$, where $\Dc$ is given in 
(\ref{Dcmat});
this implies 
$\al(T^{{\*}^{\#^{\natural_1}}})=[\al(T)]^{{\*}^{\#^{\natural_2}}}~\mbox{:}$ 
\vskip -0.3em
\[
\begin{array}{c}
{}~~~~~~\al({T^{{\*}^{\#^{\natural_1}}}})\,=\al(C^t\Dc\Dc'_1T\Dc'_1\Dc C^t)
=(C^t\Dc\Dc'_1A)T(A^{-1}\Dc'_1\Dc C^t)\nonumber\\
{}[\al(T)]^{{\*}^{\#^{\natural_2}}}=
\overline{A}\,T^{{\*}^{\#^{\natural_2}}}\,\overline{A}^{-1}=
(\overline{A}C^t\Dc\Dc'_2)T (\Dc'_2\Dc C^t\overline{A}^{-1})~.
\end{array}\nonumber
\]
{\vskip -1.8em \rightline{$_{\textstyle{\Box}}$}
\sk
We now study quantum deformations of $SO^*(2n)$.
The real form  $SO^*(2n)$ in the commutative case is generated by the matrices
$O$ that satisfy
\eq
O^tO={\bf 1} ~~~,~~~~~\overline{O}^tJO=J, ~~~~~\mbox{ where }
J=\mbox{$\scriptsize{\left(
\begin{array}{cc}
{0}&{{\bf 1}_n}\\
{-{\bf 1}_n}&{0}
\end{array}\right)}$}\label{so*}~.
\en
\vskip -0.2em
\noi Relation $\overline{O}^tJO=J$ is equivalent to $\overline{O}=JOJ^{-1}$.
There are many real forms $SO^*_q(2n)$, they are described by the matrices 
$\Dc''=i\, 
diag(\epsi_1,...\, \epsi_{2n})$ 
with $\epsi_j^2=1$,  $\epsi_{j'}=-\epsi_j$,
$j=1,\ldots , 2n$ and $\epsi_{n}=-\epsi_{n+1}=1$. 
It is easy to verify that  the matrices $\Dc''$ satisfy
(\ref{DDRDD}), (\ref{DCD}),  $\Dc''^2=-{\bf 1}$, (\ref{a}) and that therefore 
define involutive 
Hopf algebra automorphisms $\natural~:~ T\mapsto 
T^{\natural}\equiv \Dc''T\Dc''^{-1}$
that commute with $\*$ and $\xx$. The conjugations ${{\xx}^{\natural}}$ 
are equivalent to $\xx$ [hint: use 
$A=diag(e^{-i\epsi_{\!_1}\pi/_{\scriptstyle 4}},...
\,e^{-i\epsi_{\!_{2n}}\pi/_{\scriptstyle 4}})$] 
and still give the real form $SO_q(n,n)$  with $|q|=1$; 
since $\# {\scriptstyle{\, {}^{{}_{\circ}}}}\natural= 
\natural {\scriptstyle{\, {}^{{}_{\circ}}}}\#$,
the star structures ${{{\xx}^{\#}}^{\natural}}$ 
are equivalent to ${\xx}^{\#}$ and still give the real 
form $SO_q(n+1,n-1)$ with $|q|=1$.
 
The conjugations $\*^{\natural}$ give the real forms $SO_q(2n, \Dc'')$.
To prove that all the $q$-groups $SO_q(2n, \Dc'')$ are deformations of 
$SO^*(2n)$
we show  that in the commutative limit all the conjugations  $\*^{\natural}$
are equivalent. Then we find the transformation to the basis (\ref{so*}).

\noi{\sl Proof}: $\*^{\natural}$ given by 
$\Dc''=i\,diag(\epsi_1,... , \epsi_{2n})$, for $q=1$, is equivalent to the 
conjugation given by
$\Dc''_1=i\,diag(1,...1,-1,...-1)$ via 
the authomorphism $\al(T)=ATA^{-1}$ 
where $A$ is defined in (\ref{ultim}).
Notice that only in the $q=1$ limit $A$ satisfies 
(\ref{DDRDD}) and (\ref{DCD}).
\eq
{A= {\scriptsize{
{1\over 2}\left(\!  \begin{array}{cccccccc}
{\!\!1\!+\!\epsi_{\!1}}&{}&{}&{}&{}&{}&{}&{\!\!\!\!1\!-\!\epsi_{\!1}}\\
{}&{\!\!\!\!\cdot\cdot\cdot}&{}&{}&{}&{}&{\!\!\!\!\!\cdot\cdot\cdot}&{}\\
{}&{}&{\!\!\!\!1\!+\!
\epsi_{\!n\!-\!1}}&{}&{}&{\!\!\!1\!-\!\epsi_{\!n\!-\!1}}&{}&{}\\
{}&{}&{}&{\!\!\!2}&{\!\!\!0}&{}&{}&{}\\
{}&{}&{}&{\!\!\!0}&{\!\!\!2}&{}&{}&{}\\
{}&{}&{\!\!\!1\!-\!\epsi_{\!n\!-\!1}}&{}&{}&{\!\!\!1\!+
\!\epsi_{\!n\!-\!1}}&{}&{}\\
{}&{\!\!\!\cdot\cdot\cdot}&{}&{}&{}&{}&{\!\!\!\!\cdot\cdot\cdot}&{}\\
{\!\!1\!-\!\epsi_{\!1}}&{}&{}&{}&{}&{}&{}&{\!\!\!\!1\!+\!\epsi_{\!1}}
\end{array}\!\!\!\right) }}~,~~
M''=\scriptsize{
{1\over \sqrt{2}}\left(\!\!  \begin{array}{cccccccc}
{1}&{}&{}&{}&{}&{}&{}&{\!1}\\
{}&{\!\cdot\cdot\cdot}&{}&{}&{}&{}&{\!\!\cdot\cdot\cdot}&{}\\
{}&{}&{\!1}&{}&{}&{\!1}&{}&{}\\
{}&{}&{}&{\!\!1}&{\!\!\!1}&{}&{}&{}\\
{\!i}&{}&{}&{}&{}&{}&{}&{\!-i}\\
{}&{\!i}&{}&{}&{}&{}&{\!-i}&{}\\
{}&{}&{\!\!\cdot\cdot\cdot}&{}&{}&{\!\!\cdot\cdot\cdot}&{}&{}\\
{}&{}&{}&{\!\!i}&{\!\!\!-i}&{}&{}&{}
\end{array}\!\!\!\right) }}
\label{ultim}
\en
To obtain  (\ref{so*}) consider $\Dc''_1=i\,diag(1,...1,-1,...-1)$,
in the classical limit and in the basis $O= M''TM''^{-1}$, where $M''$ is 
defined in 
({\ref{ultim}}), relation $T^tCT=C$ reads
$O^tO=${\bf 1} and relation 
$T^{{\*}^\natural}=C^t\Dc''T\Dc''^{-1}C^t$ reads $\overline{O}=JOJ^{-1}$. 
[Hint: $\overline{O}=O^{{\*}^{\natural}}=\overline{M''}T^{{\*}^{\natural}}
\overline{M''}^{-1}$].$\rightline {}\!\!\!\!\!\!\!\!\!\!\!\!\!\!\!\!
\!\!\!\!\!\!\!\!\!\!\!\!\!\!\!\!\!\!\!\!\!\!\!\Box$
\sk\sk
In conclusion in the $SO_q(2n)$ case, for $q\in\Rb$, there are 
$2^{n-1}+2^{n-2}$ real forms 
$SO_q(2n, \Dc')$ and $2^{n-1}$ real forms $SO_q(2n, \Dc'')$ for a total of 
$2^{n}+2^{n-2}$
real forms as in \cite{Twietmeyer}. Also, as in \cite{Twietmeyer}, for $|q|=1$ 
we have only
$2$ inequivalent conjugations. 
\sk
We have thus contructed all real forms of $SO_q(N)$.
\sk
\noi{\bf Note 2}$~$
In the classical case 
real forms of $\mbox{\sl{so}}(N)$ (and more in general of a semisimple Lie 
algebra)
are classified via involutive automorphisms of the compact real form 
$\mbox{\sl{so}}(N,{\bf\mbox{R}})$.
(See for example \cite{class}).
The procedure we have followed to describe the real forms of $SO_q(N)$ is 
exactly the same as 
that of the classical case reformulated in the $*$-structure language.
More explicitly, we recall that an involutive automorphism $\sharp$ of the 
compact form 
$\mbox{\sl g}_{ct}$ of the semisimple complex Lie algebra 
$\mbox{\sl g}=\mbox{\sl g}_{ct}\oplus i\, \mbox{\sl g}_{ct}$,
splits $\mbox{\sl g}_{ct}$ in the direct sum $\mbox{\sl g}_{ct}=\mbox{\sl t}
\oplus\mbox{\sl p}$
where $\mbox{\sl t, p}$ are respectively the  $\sharp$ eigenspaces with 
eigenvalue $+1,-1$.
Then  $\mbox{\sl g}_{\Rbold}=\mbox{\sl t}\oplus i\,\mbox{\sl p}$, is 
a non-compact real form of $\mbox{\sl g}$.
If $\*$ is the star structure associated to $\mbox{\sl g}_{ct}\,:~ 
\mbox{\sl g}_{ct}=\{\chi\in \mbox{\sl g}~|~{\chi^{\*}}=-\chi\}$ (recall 
Section 2),
then we have $\mbox{\sl g}_{\Rbold}=\mbox{\sl t}\oplus i\,\mbox{\sl p}=
\{\chi\in \mbox{\sl g}~|~{\chi^{\*}}^{\#}=-\chi\}$;
this last relation  shows that ${\*}^{\#}$ is  the star structure $*$
canonically associated to the real form $\mbox{\sl g}_{\Rbold}$. It naturally
extends to ${\cal U}(\mbox{\sl g})$ so that 
\eq
{}~~~~~~~~~~~~~~~~~~~~~~~~~~~*=\# {\scriptstyle{\; 
{}^{{}_{\circ}}}}\!\*\label{unno}~~~~~~~~~~~~~~~~~~~~~ 
\mbox{ on }~{\cal U}(\mbox{\sl g})~.
\en
Since all real forms can be found
via involutive automorphisms of $\g_{ct}$, all $*$-structures can be found 
via (\ref{unno})
i.e. via  composition of the compact star structure $\*$ with involutive 
authomorphisms of 
$\mbox{\sl g}_{ct}$. Notice that $\#$ is an involutive automorphism of 
$\mbox{\sl g}_{ct}$
if and only if $\#$ is an involutive automorphism of 
$\mbox{\sl g}=\mbox{\sl g}_{ct}\oplus i\, \mbox{\sl g}_{ct}$ and 
$\# {\scriptstyle{\, {}^{{}_{\circ}}}}\*= 
\* {\scriptstyle{\, {}^{{}_{\circ}}}}\#$.

In the orthogonal case (for $N\not=8$) 
all automorphisms are realized, in the defining representation of 
$\mbox{\sl g}_{ct}=\mbox{\sl{so}}(N,\Rb)$, via an 
orthogonal matrix $\Dc\in O(N,\Rb)$ (with $\Dc^2=\pm \mbox{\bf 1}$ if 
the automorphism is involutive). 
Explicitly, let $\{\chi_i\}$ be a basis of $\g_{ct}$, we have 
$\#$: $\chi_i\mapsto\chi_i^{\#}\equiv \Dc\chi_i \Dc^{-1}\equiv \
{\sl\bf D}^{\;j}_{i}\chi_j$. We can use the duality (\ref{U*A}) to deduce 
how $*$ in 
(\ref{unno})
acts on the functionals $\T{a}{b}\in {\cal{F}}(G)$ [$G=SO(N,\Cb)$].  
More simply, apply (\ref{unno})
to  $\T{a}{b}(g)=({e^{x^i(g)\chi_i}})^a{}_b$ where $g\in G$ and $x^i(g)$ are the (complex) 
coordinates of $g$; we obtain
\eq
T^*=T^{\*^{\#}}=\Dc T^{\*}\Dc^{-1}\label{ddue}~.
\en
This formula  shows that in the classical case all $*$-structures of 
${\cal{F}}(G)$ are given
via involutive automorphisms of $\mbox{\sl g}_{ct}=\mbox{\sl{so}}(N,\Rb)$, 
i.e.,
via involutive automorphisms of $Fun(SO(N,\Rb))$; they are realized 
via matrices $\Dc\in O(N,\Rb)$. (For $N=8$ all $*$-structures 
are given by (\ref{ddue}) as well).
The same occurs in the quantum case. 
The $*$-structures $\*^{\#}$, $\*^{\natural}$ and $\*^{\#^{\natural}}$  
we have considered 
have all the form (\ref{ddue}) and correspond to the involutive authomorphisms 
${\#}$, ${\natural}$ and ${\#^{\natural}}$ of $SO_q(N)$ that commute 
with $\*$ (see Theorem 1),
i.e. correspond to the involutive automorphisms of the compact form  
$SO_q(N,\mbox{\bf R})$. 
${\#}$, ${\natural}$ and ${\#^{\natural}}$ are realized via the matrices 
$\Dc',\Dc'',\Dc\Dc'$ [$\Dc$ here is given in (\ref{Dcmat})]; these  
are compact quantum orthogonal matrices, i.e. are matrices with complex 
entries that satisfy the orthogonality and reality relations:
(\ref{RTT}), (\ref{Torthogonalitymat}), 
(\ref{a}) [cf. Note 1] and have $q$-determinant$\;=\pm 1$.  

{$\!\!$}We now show that $\Dc',\Dc'',\Dc\Dc'$ classify the involutive 
automorphisms of
$SO_q(N,\Rb)$ ($\forall N\not=8$).
From \cite{Twietmeyer} we deduce that $SO_q(2n+1,\Rb)$
has $2^n$ involutive automorphisms. 
The involutive automorphisms given by $T^\natural=\Dc'T\Dc'^{-1}$ 
[$\Dc'={{diag}(\epsi_1,\ldots ,\epsi_{2n+1})}$; recall that
$-\Dc$, where  $\Dc$ is defined in (\ref{Dcmat}) is a particular $\Dc'$ 
matrix]  
are $2^n$ and we conclude that 
these are all the involutive automorphisms of $SO_q(2n+1,\Rb)$.
In  the classical case
all automorphisms of $\mbox{\sl{so}}(2n+1,\Rb)$ are inner i.e. 
can be realized via 
matrices  of $SO(2n+1,\Rb)$, similarly, 
in the quantum case,
the matrices $\Dc'={{diag}(\epsi_1,\ldots ,\epsi_{2n+1})}$
satisfy the defining relations of the
compact quantum group $SO_q(2n+1,{\mbox{\bf R}})$.

The automorphism group of $\mbox{\sl{so}}(2n,\Rb)$ 
(for $n\not=4$) has two disconnected components and 
is $O(2n,\Rb)$.
This corresponds to the discrete {\bf Z}$_2$
symmetry of the Dynkin diagram of the $D_n$ series.
In the quantum case there are $2^n$ \cite{Twietmeyer} involutive inner 
automorphisms. 
These are all  
realized via the $2^{n-1}$ matrices $\Dc'$
and  the $2^{n-1}$ matrices $\Dc''$; these matrices satisfy the defining 
relations 
of  $SO_q(N,{\mbox{\bf R}})$.
On the other hand there are $2^{n-1}$ outer involutive authomorphisms,
these are all realized via the $2^{n-1}$ matrices $\Dc\Dc'$
that are quantum orthogonal matrices with
 det$_q\,\Dc\Dc'=-1$.
\sk

\noi{\bf Note 3}$~$
There are two different Lorentz groups obtained as real forms of $SO_q(4)$:
$SO_q(3,1)$ with $q\in {\bf\mbox{R}}$ and $SO_q(3,1)$ with $|q|=1$.
More in general we have deformations of  $SO(n+1,n)$, $SO(n,n)$, $SO(n+1,n-1)$ 
both for $|q|=1$ and for $q\in {\bf\mbox{R}}$. 
\sk
\noi{\bf Note 4}$~$
As shown in \cite{altroarticolo} 
the classical embedding of $SO(N)$ in $SO(N+2)$ holds also in the quantum case.
We have that the universal enveloping algebra ${\cal U}_q(so(N))$ of $SO_q(N)$
is a Hopf subalgebra of ${\cal U}_q(so(N+2))$. Otherwise stated 
$SO_q(N)$, that is dual to  ${\cal U}_q(so(N))$, is a quotient  of $SO(N+2)$
via a Hopf ideal $H$:
$SO_q(N)=SO(N+2)/H$. In the commutative case $H$ is the ideal in $Fun(SO(N+2))$ 
generated by the equivalence relation $f\sim f'\Leftrightarrow 
f|_{SO(N)}=f'|_{SO(N)}$ where
$f$ and $f'$ are generic functions on $SO(N+2)$ and $f|_{SO(N)}$ is $f$ 
restricted
to ${SO(N)}$.

On the other hand the classical embedding of $SO(N)$ in $SO(N+1)$ does 
not occur at the 
quantum level. 
In particular the $q$-Lorentz groups discussed here
do not ``contain" three dimensional $q$-orthogonal subgroups, i.e., $SO_q(3)$ 
cannot be found as a quotient $SO_q(4)/H$ with respect to an appropriate 
ideal $H$.  

It is interesting to note that however the $q$-Minkowski plane for $q\in\Rb$ 
``contains" the euclidean $q$-plane and the three dimensional $q$-Minkowski 
plane.

\noi {\sl Proof} :$~$ The $q$-Minkowski coordinates satisfy:
\eq
\begin{array}{c}
x^1x^2=qx^2x^1~,~
x^1x^3=qx^3x^1~,~
x^2x^4=qx^4x^2~,~
x^3x^4=qx^4x^3~,~\\
x^2x^3=x^3x^2~,~~x^1x^4=x^4x^1-(q-q^{-1})x^2x^3~.
\end{array}\label{4plane}
\en
with $(x^1)^{\b*}=q x^4$, $(x^2)^{\b*}=x^2$ and  $(x^3)^{\b*}=x^3$.
The $SO_q(3)$ plane commutations are
\eq
y^1y^2=qy^2y^1~,~
y^2y^3=qy^3y^2~,~y^1y^3=y^3y^1-(q^{1\over 2}-q^{-{1\over 2}})y^2y^2
\label{3plane}
\en
The conjugation that gives the euclidean plane is 
$(y^1)^{\*}=qy^3$, $(y^2)^{\*}=y^2$, while for the $SO_q(2,1)$ plane 
we have $(y^1)^{\b*}=qy^3$, $(y^2)^{\b*}=-y^2$.
To obtain the $SO_q(3,{\bf\mbox{R}})$ plane from (\ref{4plane}) we quotient 
with respect to the relation $x^2=x^3$ and identify 
$x^1,\sqrt{q^{1\over 2}+q^{-{1\over 2}}}x^2,x^4$ respectively with $y^1,y^2,y^3$.
To obtain the $SO_q(2,1)$ plane we impose the relation $x^2=-x^3$ 
and identify $x^1,i\sqrt{q^{1\over 2}+q^{-{1\over 2}}}x^2,x^4$ respectively 
with $y^1,y^2,y^3\,$. \QED

One similarly proves that the $|q|=1$ $SO_q(2,1)$ plane is obtained 
from the $|q|=1$ $q$-Minkowski spacetime with 
${{\xx}^{\sharp}}$-conjugation (\ref{plane}): $({x^1}){{\x}^{{}^{\!\!\!\sharp}}}=x^1, 
({x^2}){{\x}^{{}^{\!\!\!\sharp}}}=x^3, ({x^4}){{\x}^{{}^{\!\!\!\sharp}}}=x^4$
(cf. \cite{Firenze1}), via the quotient $x^2=x^3$. Note that this $|q|=1$ 
$q$-Minkowski spacetime has also been studied in \cite{d}, there the $q$-commutations 
and the conjugation are derived from the conformal group $SU_q(2,2)$. 
Explicitly, in \cite{d} the conjugation used is ${{\xx}^{\sharp^\natural}}$ 
with ${\cal{D'}}=diag(-1,1,1-1)$:
$({x^1}){{\x}^{{}^{\!\!\!\sharp^\natural}}}=-x^1, 
({x^2}){{\x}^{{}^{\!\!\!\sharp^\natural}}}=x^3, 
({x^4}){{\x}^{{}^{\!\!\!\sharp^\natural}}}=-x^4$, this conjugation 
we have proved to be equivalent to ${{\xx}^{\sharp}}$.
\sk
\vskip .9em
{\bf {Acknowledgments}}
\sk
The autor thanks  Leonardo Castellani for fruitful discussions and 
encouragement. 
The autor thanks Bogdan Morariu, Nicolai Reshetikhin and Bruno Zumino for their
valuable 
comments and suggestions.
\sk
This work has been supported  by an INFN post-doctoral 
grant (concorso n. 6077/96). 
It has been accomplished through the Director,
Office of Energy Research, Office af High Energy and Nuclear  
Physics, Division of High Energy Physics
of the U.S. Department of Energy under Contract DE-AC03-76SF00098
and by the National Science Foundation under grant PHY-95-14797.

\vfill\eject

\begin{thebibliography}{99}

\bibitem{Wess} B. L. Cerchiai, J. Wess,   {\sl q-Deformed Minkowski Space 
based on a 
q-Lorentz Algebra}
math.QA/9801104;
J. Wess, {\sl Quantum groups and $q$-lattices in phase-space}, in the
Proceedings of the 5-th Hellenic School and Workshop on Elementary
Particle Physics, Corfu, September 1995, q-alg/9607002.

\bibitem{FRT} L.D. Faddeev,
N.Yu. Reshetikhin and L.A. Takhtajan, {\sl Quantization of Lie Groups and 
Lie Algebras},
Algebra i Anal. 1 {\bf 1} (1989) 178 (Leningrad Math. J. {\bf 1} 193
(1990)). 

\bibitem{Firenze1} E. Celeghini, R. Giachetti, A. Reyman, E. Sorace,
 M. Tarlini, {\sl $SO_q(n+1,N-1)$ as a Real Form of $SO_q(2n,\Cb)$},
 Lett. Math. Phys. {\bf 23} (1991) 45.

\bibitem{Dobrev}
V.K. Dobrev, {\sl Canonical $q$-Deformations of 
Noncompact Lie (Super-) Algebras}, J. Phys. A: Math. Gen. {\bf 26} (1993)
1317, (G\"{o}ttingen University preprint, July 1991).  

\bibitem{Twietmeyer} E. Twietmeyer, {\sl Real Forms of 
${\cal{U}}_q(\mbox{\sl{g}})$},
 Lett. Math. Phys. {\bf 24} (1992) 49.

\bibitem{Fiore} G. Fiore, J. Phys. {\sl Quantum Groups $SO_q(N)$, $Sp_q(N)$
have determinants, too}
{\bf A27} (1994) 3795.

\bibitem{d}  V.K. Dobrev, {\sl New $q$ - Minkowski space-time and $q$ - Maxwell
 equations hierarchy from $q$-conformal invariance} Phys. Lett.
 {\bf 341B} 133 (1994), erratum -ibid. {\bf 346B} 427 (1995).                 

\bibitem{altroarticolo} P. Aschieri and  L. Castellani,
{\sl Universal enveloping algebra and differential calculi on 
inhomogeneous orthogonal $q$-groups}, LBNL 40330, q-alg/9705023.
to be publ. in  J. Geom. and Phys.

\bibitem{class} R. Gilmore {\sl Lie Groups, Lie algebras and Some of Their 
Applications} 
Wiley-Interscience  New York - London - Sydney - Toronto (1974).  M. Hausner, 
J. Schwartz, 
{\sl Lie Groups; Lie Algebras}
Gordon and Breach, New York - London - Paris (1968).
\end{thebibliography}
\end{document}